\documentclass[draft,numbook]{svjour3}

\usepackage{graphicx}
\usepackage{url}
\usepackage{amssymb}
\smartqed  

\begin{document}

\title{On the Hurwitz Zeta Function}

\author{Lazhar Fekih-Ahmed}

\institute{Lazhar Fekih-Ahmed \at
         \'{E}cole Nationale d'Ing\'{e}nieurs de Tunis, BP 37, Le Belv\'{e}d\`{e}re 1002 , Tunis, Tunisia \\
              \email{lazhar.fekihahmed@enit.rnu.tn}           
}

\date{May 2, 2012
 }

\maketitle

\begin{abstract}
We give new integral and series representations of the Hurwitz
zeta function. We also provide a closed-form expression of the
coefficients of the Laurent expansion of the Hurwitz-zeta function
about any point in the complex plane.

\end{abstract}

\section{Introduction}
\label{intro}

The Hurwitz zeta function defined by the series
\begin{equation}\label{sec1-eq1} \zeta(s,a)=
\frac{1}{a^{s}}+\frac{1}{(1+a)^{s}}+\cdots+\frac{1}{(k-1+a)^{s}}+\cdots=\sum_{n=1}^{\infty}(n-1+a)^{-s},
\end{equation}
where $0< a\le 1$, is a well-defined series when $\Re(s)>1$, and
can be analytically continued to the whole complex plane with one
singularity, a simple pole with residue $1$ at $s=1$.

The Hurwitz zeta function has also the following integral
representation

\begin{equation}\label{sec1-eq2}
\zeta(s,a)=\frac{1}{\Gamma(s)}\int_{0}^{\infty}
\frac{e^{-(a-1)t}}{e^{t}-1} t^{s-1}\,dt,
\end{equation}

valid for $\Re(s)>1$ \cite{whittaker}. Moreover, it has the
following analytic continuation represented by the following
contour integral
\begin{equation}\label{sec1-eq3}
\frac{\Gamma(1-s)}{2\pi
i}\int_{\mathcal{C}}\frac{e^{(a-1)t}}{e^{-t}-1} t^{s-1}\,dt,
\end{equation}

where $\mathcal{C}$ is the Hankel contour consisting of the three
parts $C=C_{-}\cup C_{\epsilon}\cup C_{+}$:  a path  which extends
from $(-\infty,-\epsilon)$, around the origin counter clockwise on
a circle of center the origin and of radius $\epsilon$ and back to
$(-\epsilon,-\infty)$, where $\epsilon$ is an arbitrarily small
positive number.

The integral (\ref{sec1-eq3})  defines $\zeta(s,a)$ for all
$s\in\mathbb{C}$ with a single pole at $s=1$.

\section{New Integral and Series representations of $\zeta(s,a)$}\label{sec2}

The new representations of $\zeta(s,a)$ are based on the
alternating sums defined by

\begin{equation}\label{sec2-eq1}
S_{n}(s,a) = \sum_{k=0}^{n-1}(-1)^{k}{n-1\choose k}(k+a)^{-s}~{\rm
for}~ n\ge 2,
\end{equation}

and by $S_{1}(s,a)=1$ for $n=1$, and on the real function
$\psi(t)$ defined by

\begin{equation}\label{sec2-eq1bis}
\psi(t)=\frac{te^{t}}{(e^{t}-1)^2}-\frac{1}{e^{t}-1}+\frac{(a-1)t}{e^{t}-1}.
\end{equation}
We prove the following theorem

\begin{theorem}\label{sec2-thm1}
For all $s$ such that $\Re(s)>0$ and all $0< a\le1$, we have

\begin{eqnarray*}
  \rm{(A)} \qquad\qquad&& (s-1)\zeta(s,a)
=\frac{1}{\Gamma(s)}\int_{0}^{\infty}\psi(t) e^{-(a-1)t}t^{s-1}\,dt.\\
  \rm{(B)} \qquad\qquad && (s-1)\zeta(s,a)
=\sum_{n=1}^{\infty}S_{n}(s,a)\bigg(\frac{1}{n+1}+
\frac{a-1}{n}\bigg).
\end{eqnarray*}

\end{theorem}
\begin{proof}

We can rewrite $S_{n}(s,a)$ in (\ref{sec2-eq1}) as

\begin{eqnarray}
S_{n}(s,a) &=&
\frac{1}{\Gamma(s)}\int_{0}^{\infty}\sum_{k=0}^{n-1}(-1)^{k}{n-1\choose
k}
   e^{-(k+a)t}t^{s-1}\,dt\nonumber \\
  &=&
  \frac{1}{\Gamma(s)}\int_{0}^{\infty}(1-e^{-t})^{n-1}e^{-at}t^{s-1}\,dt,\label{sec2-eq3}
\end{eqnarray}

since we know that

\begin{equation}\label{sec2-eq2}
(n+a-1)^{-s}=\frac{1}{\Gamma(s)}\int_{0}^{\infty}e^{-(n+a-1)t}t^{s-1}\,dt,
\end{equation}

a valid formula for $\displaystyle{\Re(s)>0}$, and since

\begin{equation}\label{sec2-eq2bis}
\sum_{k=0}^{n-1}(-1)^{k}{n-1\choose k}
e^{-(k+a)t}=e^{-at}(1-e^{-t})^{n-1}.
\end{equation}

The proof consists of an evaluation of the  sum

\begin{equation}\label{sec2-eq4}
\sum_{n=1}^{\infty}\frac{S_{n}(s,a)}{n+1}+(a-1)\sum_{n=1}^{\infty}\frac{S_{n}(s,a)}{n}.
\end{equation}

But before we do so, we need to establish the following two
identities valid for $0<t<\infty$:

\begin{eqnarray}
\frac{te^{-t}}{1-e^{-t}}&=&\sum_{n=1}^{\infty}\frac{(1-e^{-t})^{n-1}e^{-t}}{n}.\label{sec2-eq8}\\
\frac{te^{-t}}{(1-e^{-t})^2}-\frac{e^{-t}}{1-e^{-t}}&=&\sum_{n=1}^{\infty}\frac{(1-e^{-t})^{n-1}e^{-t}}{n+1}.
\label{sec2-eq9}
\end{eqnarray}

To prove these identities, we start from the series

\begin{equation}\label{sec2-eq5}
t=-\log(1-(1-e^{-t}))=\sum_{n=1}^{\infty}\frac{(1-e^{-t})^{n}}{n}
\end{equation}
which is valid for $0<t<\infty$.

This yields

\begin{equation}\label{sec2-eq6}
\frac{t}{1-e^{-t}}=\sum_{n=1}^{\infty}\frac{(1-e^{-t})^{n-1}}{n},
\end{equation}

and

\begin{equation}\label{sec2-eq7}
\frac{t}{(1-e^{-t})^2}=\frac{1}{1-e^{-t}}+\sum_{n=1}^{\infty}\frac{(1-e^{-t})^{n-1}}{n+1},
\end{equation}

from which the two identities follow easily.

Without worrying about interchanging sums and integrals for the
moment, we have

\begin{eqnarray}\label{sec2-eq10}
\sum_{n=1}^{\infty}\frac{S_{n}(s,a)}{n+1}&=&\frac{1}{\Gamma(s)}\sum_{n=1}^{\infty}\int_{0}^{\infty}\frac{(1-e^{-t})^{n-1}}{n+1}e^{-at}t^{s-1}\,dt\\
&=&\frac{1}{\Gamma(s)}\int_{0}^{\infty}\sum_{n=1}^{\infty}\frac{(1-e^{-t})^{n-1}}{n+1}e^{-at}t^{s-1}\,dt\label{sec2-eq11}\\
&=& \frac{1}{\Gamma(s)}\int_{0}^{\infty}\bigg
(\frac{te^{-t}}{(1-e^{-t})^2}-\frac{e^{-t}}{1-e^{-t}}\bigg )
e^{-(a-1)t}t^{s-1}\,dt,\label{sec2-eq12}
\end{eqnarray}

and similarly

\begin{eqnarray}\label{sec2-eq13}
\sum_{n=1}^{\infty}\frac{S_{n}(s,a)}{n}&=&\frac{1}{\Gamma(s)}\sum_{n=1}^{\infty}\int_{0}^{\infty}\frac{(1-e^{-t})^{n-1}}{n}e^{-at}t^{s-1}\,dt\\
&=&\frac{1}{\Gamma(s)}\int_{0}^{\infty}\sum_{n=1}^{\infty}\frac{(1-e^{-t})^{n-1}}{n}e^{-at}t^{s-1}\,dt\label{sec2-eq14}\\
&=& \frac{1}{\Gamma(s)}\int_{0}^{\infty}
\frac{te^{-t}}{1-e^{-t}}e^{-(a-1)t}t^{s-1}\,dt.\label{sec2-eq15}
\end{eqnarray}

Therefore,

\begin{eqnarray}
&&\sum_{n=1}^{\infty}\frac{S_{n}(s,a)}{n+1}+(a-1)\sum_{n=1}^{\infty}\frac{S_{n}(s,a)}{n}=\nonumber\\
&&\frac{1}{\Gamma(s)}\int_{0}^{\infty}\bigg
(\frac{te^{t}}{(e^{t}-1)^2}-\frac{1}{e^{t}-1}+\frac{(a-1)t}{e^{t}-1}\bigg
) e^{-(a-1)t}t^{s-1}\,dt.\label{sec2-eq16}
\end{eqnarray}

This proves that the right hand sides of (A) and (B) in the
statement of the theorem are equal. Now, by observing that
\begin{equation}\label{sec2-eq17}
\frac{d}{dt}\bigg (\frac{-te^{-(a-1)t}}{e^{t}-1}\bigg )=\bigg
(\frac{te^{t}}{(e^{t}-1)^2}-\frac{1}{e^{t}-1}+\frac{(a-1)t}{e^{t}-1}\bigg
) e^{-(a-1)t},
\end{equation}

we can perform an  integration by parts in (\ref{sec2-eq16}) when
$\Re(s)>1$. The integral in the right hand side of
(\ref{sec2-eq16}) becomes

\begin{equation}\label{sec2-eq18}
\frac{s-1}{\Gamma(s)}\int_{0}^{\infty} \frac{e^{-(a-1)t}}{e^{t}-1}
t^{s-1}\,dt=(s-1)\zeta(s,a).
\end{equation}

Thus,

\begin{equation}\label{sec2-eq19}
\sum_{n=1}^{\infty}\frac{S_{n}(s,a)}{n+1}+(a-1)\sum_{n=1}^{\infty}\frac{S_{n}(s,a)}{n}=(s-1)\zeta(s,a),
\end{equation}

and this prove the theorem when $\Re(s)>1$. However, formula
(\ref{sec2-eq19}) remains valid for $\Re(s)>0$ since the integral
(\ref{sec2-eq16}) is well-defined for $\Re(s)>0$.

To finish the proof, we now need to justify the interchange of
summation and integration in both equation (\ref{sec2-eq11}) and
equation (\ref{sec2-eq14}). We show this justification for
equation (\ref{sec2-eq11}) only, the other is similar.

The interchange in (\ref{sec2-eq11}) is indeed valid because the
series
\begin{equation}\label{sec2-eq20}
\sum_{n=1}^{\infty}\int_{0}^{\infty}\frac{(1-e^{-t})^{n-1}}{n+1}e^{-at}t^{s-1}\,dt
\end{equation}

converges absolutely and uniformly for $0<t<\infty$. To prove
this, it suffices to  show uniform convergence for the dominating
series
\begin{equation}\label{sec2-eq21}
\sum_{n=1}^{\infty}\int_{0}^{\infty}\frac{(1-e^{-t})^{n-1}}{n+1}e^{-at}t^{\sigma-1}\,dt,
\end{equation}
where $\sigma=\Re(s)$.

Indeed, let $K={\rm max}((1-e^{-t})^{n-1}e^{-t/2})$, $0<t<\infty$.
A straightforward calculation of the derivative shows that
\begin{equation}\label{sec2-eq22}
K=(1-\frac{1}{2n-1})^{n-1}\frac{1}{\sqrt{2n-1}}
\end{equation}
 and  is attained
when $e^{-t}=\frac{1}{2n-1}$. Now, for $n\ge 2$, we have
\begin{eqnarray}
\int_{0}^{\infty}(1-e^{-t})^{n-1}e^{-at}t^{\sigma-1}\,dt&=&
\int_{0}^{\infty}(1-e^{-t})^{n-1}e^{-t/2}
(e^{-(-1/2+a)t}t^{\sigma-1})\,dt\nonumber\\
&\le &K \int_{0}^{\infty}e^{-(-1/2+a)t}t^{\sigma-1}\,dt\nonumber\\
&=&(1-\frac{1}{2n-1})^{n-1}\frac{\Gamma(\sigma)}{\sqrt{2n-1}(-1/2+a)^{\sigma}}\nonumber\\
&\le&\frac{K^{\prime}}{\sqrt{2n-1}}.\label{sec2-eq23}
\end{eqnarray}

The last inequality implies that each term of the dominating
series is bounded by $K^{\prime}/(n+1) \sqrt{2n-1}$. Thus the
dominating series converges by the comparison test. This completes
the proof of the theorem. \smartqed \qed
\end{proof}

With $\mathcal{C}$ being the Hankel contour defined previously, an
immediate corollary of Theorem~\ref{sec2-thm1} is the following
\begin{corollary}\label{sec2-coro1}
For $0< a\le1$ and for all $s\in \mathbb{C}$, we have

\begin{eqnarray*}
  \rm{(C)} \qquad\qquad&& (s-1)\zeta(s,a)
=\frac{\Gamma(1-s)}{2\pi
i}\int_{\mathcal{C}}\psi(-t)e^{(a-1)t}t^{s-1}\,dt.\\
  \rm{(D)} \qquad\qquad && (s-1)\zeta(s,a)
=\sum_{n=1}^{\infty}S_{n}(s,a)\bigg(\frac{1}{n+1}+
\frac{a-1}{n}\bigg).
\end{eqnarray*}
\end{corollary}
\begin{proof}
The proof the first statement follows the same steps as in
\cite{whittaker} for example. As for the second statement, we can
either proceed as in \cite{hasse} or as in
\cite{lazhar}\footnote{The notation used in \cite{lazhar} may
cause some confusion. Contrary to the notation of this paper in
which $\zeta(s,a)$ denotes the Hurwitz zeta function, $\zeta(s,x)$
in the cited paper is a power series associated with $\zeta(s)$
and has nothing to do with the Hurwitz zeta function.}. In
\cite{lazhar}, the proof was given for $\zeta(s)$ (i.e. $a=1$) and
uses an estimate of the exact asymptotic order of growth of
$S_{n}(s,1)$ when $n$ is large. By looking at the definition of
$S_{n}(s,a)$, we can see that $S_{n}(s,1)$ are the Stirling
numbers of the second kind modulo a multiplicative factor when
$s\in \{0,-1,-2,\cdots\}$. Therefore, when $s=-k$, $k$ a positive
integer, $S_{n}(-k,1)$ are eventually zero for $n$ large enough.
Similarly, $S_{n}(s,a)$ are the generalized Stirling numbers of
the second kind whose generating function is given by
\cite{carlitz}:

\begin{equation}\label{sec2-eq23bis}
    \sum_{k=0}^{\infty}S_{n}(-k,a)\frac{t^k}{k!}=e^{a t}(1-e^{t})^{n}.
\end{equation}

It follows that

\begin{equation}\label{sec2-eq23bisbis}
   S_{n}(-k,a)=\frac{d^{k}}{dt^{k}}\bigg\{ e^{a
   t}(1-e^{t})^{n}\bigg\}\bigg |_{t=0},
\end{equation}

and hence $S_{n}(-k,a)$ are eventually zero for $n$ large enough.

For $s\notin \{0,-1,-2,\cdots\}$, an asymptotic estimate of
$S_{n}(s,a)$ can be obtained for $n$ large. Indeed, we have by
definition

\begin{equation}
S_{n}(s,a)=\sum_{k=0}^{n-1}{n-1\choose
k}(-1)^{k}(k+a)^{-s}\label{sec2-eq24}.
\end{equation}

By making the change of variable, $k=m-1$, the sum can  obviously
be put into  the form

\begin{equation}
S_{n}(s,a) = \frac{-1}{n}\sum_{m=1}^{n}{n\choose
m}(-1)^{m}\frac{m}{(m+a-1)^s}.\label{sec2-eq25}
\end{equation}

Asymptotic expansions of sums of the form (\ref{sec2-eq25}) have
been thoroughly studied in \cite{flajolet}. The function
$\frac{z}{(z+a-1)^s}$ has a non-integral algebraic singularity at
$s_0=1-a$ since $0<a<1$; thus, when $s$ is nonintegral,
$S_{n}(s,a)$ has the following asymptotics when $n$ is large

\begin{equation}\label{sec2-eq26}
S_{n}(s,a)\thicksim
\frac{\Gamma(1-a)n^{1-a}(\log{n})^{s-1}}{n\Gamma(s)}=\frac{\Gamma(1-a)}{n^{a}(\log{n})^{1-s}\Gamma(s)},
\end{equation}

and when $s=k\in \{1,2,\cdots\}$ the following expansion

\begin{equation}\label{sec2-eq26bis}
S_{n}(k,a)\thicksim
\frac{\Gamma(1-a)(\log{n})^{k-1}}{n^{a}(k-1)!}.
\end{equation}

For the case $a=1$ and $s$ nonintegral, the numerator becomes
equal to 1 and we obtain

\begin{equation}\label{sec2-eq27}
S_{n}(s) \thicksim \frac{1}{n(\log n)^{1-s}\Gamma(s)}.
\end{equation}

The asymptotic estimates  (\ref{sec2-eq26}), (\ref{sec2-eq26bis})
and (\ref{sec2-eq27}) are valid for $n$ large enough and for all
$s$ such that $\Re(s)>0$. To finish the proof of the corollary, we
note that the logarithmic test of series implies that our series
is dominated by a uniformly convergent series for all finite $s$
such that $\Re(s)>0$. Now, by Weierstrass theorem of the
uniqueness of analytic continuation, the function
$(s-1)\zeta(s,a)$ can be extended outside of the domain $\Re(s)>1$
and that it does not have any singularity when $\Re(s)>0$.
Moreover, by repeating the same process for $\Re(s)>-k$, $k\in
\mathbb{N}$, it is clear that the series defines an analytic
continuation of $\zeta(s)$ valid for all $s\in\mathbb{C}$.

\smartqed \qed
\end{proof}

We can also obtain the analytic continuation of $\zeta(s,a)$ to
the whole complex plane via the following

\begin{corollary}\label{sec2-coro2}
For all $s$ such that $\Re(s)>-k$ and all $0< a\le1$, we have
\begin{equation}\label{sec2-eq28}
(s-1)\zeta(s,a)
=\frac{(-1)^k}{\Gamma(s+k)}\int_{0}^{\infty}\frac{d^{k}}{dt^k}(\psi(t)
e^{-(a-1)t})t^{s+k-1}\,dt.
\end{equation}

\end{corollary}
\begin{proof}
\smartqed
 For all $k$ and all $0< a\le1$, the integrand is bounded
and has finite values at the limits of integration. Repeated
integration by parts proves the corollary. \qed
\end{proof}

Now if we consider the function $\zeta(s,a)-a^{-s}$ dealt with in
\cite{berndt}, we obtain the following corollary

\begin{corollary}\label{sec2-coro3}
Let
$\eta(t)=\frac{te^{t}}{(e^{t}-1)^2}-\frac{1}{e^{t}-1}-\frac{at}{e^{t}-1}$.
Then, for all $s$ such that $\Re(s)>0$ and all $0\le a\le1$, we
have

\begin{eqnarray*}
  \rm{(E)} \qquad\qquad&& (s-1)(\zeta(s,a)-a^{-s})
=\frac{1}{\Gamma(s)}\int_{0}^{\infty}\eta(t) e^{-at}t^{s-1}\,dt.\\
  \rm{(F)} \qquad\qquad && (s-1)(\zeta(s,a)-a^{-s})
=\sum_{n=1}^{\infty}S_{n}(s,a)\bigg(\frac{1}{n+1}-
\frac{a}{n}\bigg).
\end{eqnarray*}
\end{corollary}
\begin{proof}
\smartqed
 The proof follows the same lines as the proof of Theorem~\ref{sec2-thm1}. \qed
\end{proof}

We note that when $a=0$, $\zeta(s,a)-a^{-s}=\zeta(s)$. The
formulas of Corollary~\ref{sec2-coro3} reduce to the formulas
found in \cite{lazhar}.

\section{A formula for the Laurent coefficients}\label{sec3}

Berndt \cite{berndt} (see also the references therein) derived
expressions for the coefficients of the Laurent expansion of the
Hurwitz zeta-function $\zeta(s,a)$ about $s=1$. He also provided a
method to calculate and estimate these coefficients. In this
section, we give exact estimates of the coefficients of
$(s-1)\zeta(s,a)\Gamma(s)$ for any point of the complex plane, and
from these coefficients the Laurent coefficients of $\zeta(s,a)$
about the same point can be easily calculated.

We will only consider expansions around a point $s_0=x+iy$ in the
right half plane, i.e. $x>0$. For the other points, the reader
will realize that the extension can be easily accomplished.

The idea is to simply  use the integral formula given in formula
(A) of Theorem~\ref{sec2-thm1}. From this formula, the analytic
function $(s-1)\zeta(s)\Gamma(s)$ can be represented by a Taylor
series around any point $s_0$

\begin{equation}\label{sec3-eq1}
(s-1)\zeta(s,a)\Gamma(s)=\sum_{n=0}^{\infty}a_{n}(s-s_0)^{n}
\end{equation}

where the coefficients $a_0=(s_0-1)\zeta(s_0,a)\Gamma(s_0)$ and
$a_n$ are given by

\begin{eqnarray}
a_n&=&\frac{1}{n!}\lim_{s\to s_0}\frac{d^n}{d s^n}\bigg \{\int_{0}^{\infty}\psi(t) e^{-(a-1)t}t^{s-1}\,dt \bigg \}\nonumber\\
&=&\frac{1}{n!}\int_{0}^{\infty}\psi(t)e^{-(a-1)t}
(\log{t})^{n}t^{s_0-1}\,dt,\label{sec3-eq2}
\end{eqnarray}

with $\psi(t)$ being given in equation (\ref{sec2-eq1bis}).

Thus, following \cite{berndt}\footnote{Actually the coefficients
$\gamma_{n}(a)$ when $s_0=1$ differ slightly from ours since in
\cite{berndt} they were defined by
$(s-1)\zeta(s,a)=\sum_{n=0}^{\infty}\gamma_{n}(a)(s-1)^{n+1}$.},
if we set

\begin{equation}\label{sec3-eq3}
(s-1)\zeta(s,a)=\sum_{n=0}^{\infty}\gamma_{n}(a,s_0)(s-s_0)^{n},
\end{equation}

and if we expand $\Gamma(s)$ around the point $s_0$ into a Taylor
series of the form

\begin{equation}\label{sec3-eq4}
\Gamma(s)=\sum_{n=0}^{\infty}c_{n}(s-s_0)^{n},
\end{equation}

where $c_0=\Gamma(s_0)$ and
\begin{eqnarray}
c_n&=&\frac{1}{n!}\lim_{s\to s_0}\frac{d^n}{d s^n}\bigg \{\int_{0}^{\infty}e^{-t}t^{s-1}\,dt \bigg \}\nonumber\\
&=&\frac{1}{n!}\int_{0}^{\infty}e^{-t}
(\log{t})^{n}t^{s_0-1}\,dt,\label{sec3-eq5}
\end{eqnarray}

we easily obtain

\begin{equation}\label{sec3-eq6}
\gamma_n(a,s_0)=-\frac{1}{c_0}\sum_{k=1}^{n}\gamma_{n-k}(a,s_0)c_{k}+a_n.
\end{equation}

To get the last formula, we have merely used the coefficient
formula of division of power series. The formula is recursive and
can be used to calculate $\gamma_n(a,s_0)$ for any $0< a\le 1$ and
any $s_0$ in the right half plane. If $s_0$ is not in the right
half plane, the integral formula (\ref{sec2-eq28}) of
Corollary~\ref{sec2-coro2} can be used instead.

\end{document}